\documentclass[11pt,a4paper]{article}
\usepackage{mathrsfs}
\usepackage{cite}
\usepackage{subfigure}
\usepackage{graphicx}
\DeclareFontFamily{U}{wncy}{}
\DeclareFontShape{U}{wncy}{m}{n}{<->wncyr10}{}
\DeclareSymbolFont{mcy}{U}{wncy}{m}{n}
\DeclareMathSymbol{\Sh}{\mathord}{mcy}{"58}

\usepackage[cmex10]{amsmath}
\usepackage{amssymb}  
\usepackage{algorithmic}

\begin{document}
\title{
A Novel Fibonacci Pattern in Pascal's Triangle
}

\author{Bernhard A. Moser\\
SCCH, Austria\\
Email: bernhard.moser@scch.at }

\maketitle
\noindent
\begin{center}{\bf Abstract}\end{center}
The Fibonacci sequence is obtained as weighted sum along the rows in the Pascal triangle by choosing
a periodic up-and-down pattern of weights from the set $\{-1,-\frac{1}{2},0, \frac{1}{2}, 1\}$.
A graphical illustration of this identity shows a novel "`beautiful"' Fibonacci pattern.

\section{Motivation} 
\begin{figure}
  \begin{center}
	      \includegraphics[width=0.60 \columnwidth]{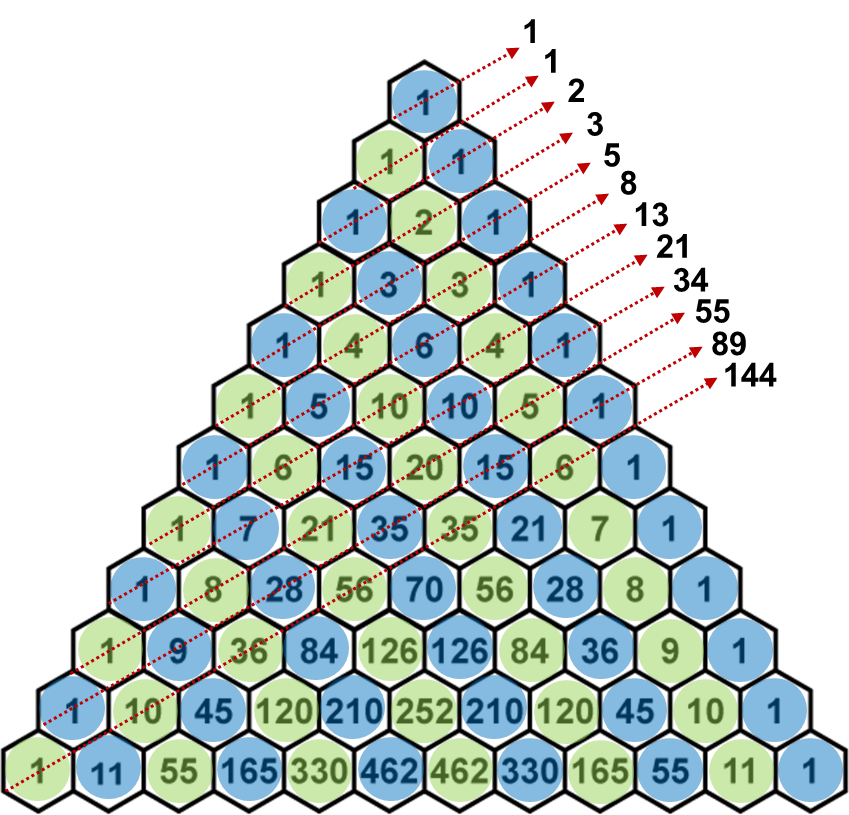} 
  \end{center}
  \caption{The classical Fibonacci pattern in the Pascal triangle.}
  \label{fig:FiboPascalKnown}
\end{figure}

Pascal's triangle and the Fibonacci numbers hide interesting and beauty patterns whose discovery has a long history 
in mathematics~\cite{bac2011,Scott2014}. 
Traditionally, these patterns play an important role in illustrating the beauty of mathematics in the classroom~\cite{Debnath11}.
One such surprising pattern is the relation between both schemes of numbers as depicted in 
Fig.~\ref{fig:FiboPascalKnown} showing that the summation along shallow diagonals in the Pascal triangle yields the Fibonacci sequence.

\begin{figure}
  \begin{center}
	      \includegraphics[width=0.90 \columnwidth]{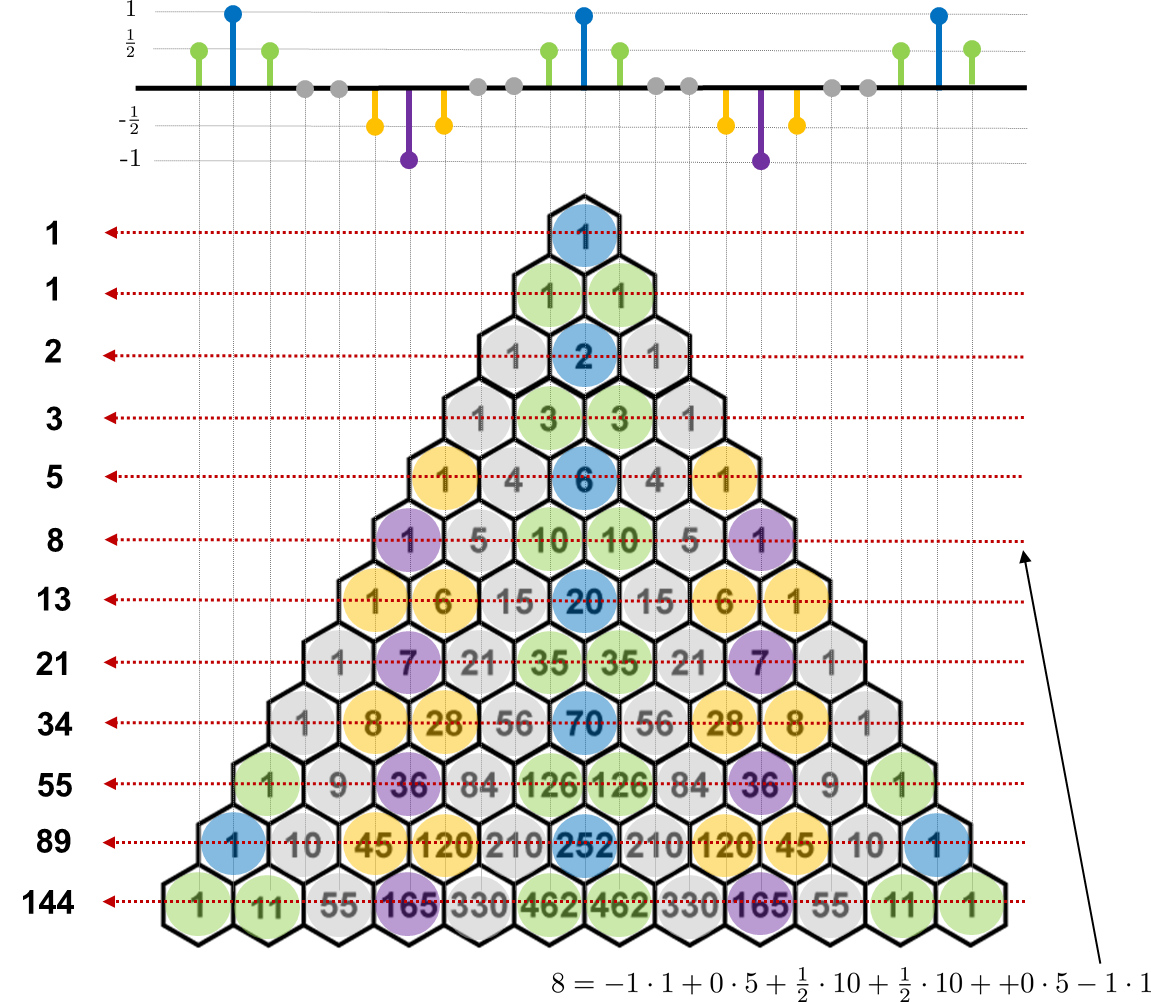} 
  \end{center}
  \caption{The novel Fibonacci pattern in the Pascal triangle.}
  \label{fig:FiboPascal}
\end{figure}

This article takes up a result from the author~\cite{Moser2014} which was a byproduct of a lattice enumeration approach.
While the identity in~\cite{Moser2014} is not at all handy, its visualization in the Pascal triangle reveals an appealing pattern as illustrated in Fig.~\ref{fig:FiboPascal}. To this end, this article is a contribution to the beauty and "`magic"' of patterns in the Pascal triangle related to the Fibonacci numbers.

\section{Sketch of Proof} 
The proof can be found in~\cite{Moser2014}.

Let 
\[ F_{k+1} = F_{k} + F_{k-1}, k\geq 2, F_1 = 1,\]
denote the Fibonacci sequence.
Further let us denote by $e_k$  the canonical unit vector in the $k$-th dimension and 
$1^T = (1,1,1,1)$.
 
The idea for the novel identity relies on the fact that
\begin{equation}
\label{eq:F}
F_{k+1} =1^T Q^k e_1,
\end{equation}
where
\begin{equation}
\label{eq:MQ}
\mathbf  Q =
\left( 								
\begin{array}{cccc}
0 & 1 & 0 & 0 \\
1 & 0 & 1 & 0 \\
0 & 1 & 0 & 1 \\
0 & 0 & 1 & 0 
\end{array}
\right). 
\end{equation}
To check Equ. (\ref{eq:F}), first observe that $R Q = Q R$ and $R = R^{-1}$,
where
\begin{equation}
\label{eq:R}
\mathbf  R =
\left( 								
\begin{array}{cccc}
0 & 0 & 0 & 1 \\
0 & 0 & 1 & 0 \\
0 & 1 & 0 & 0 \\
1 & 0 & 0 & 0 
\end{array}
\right).
\end{equation}
This means, reversing the ordering of the rows of $Q$ has the same effect as reversing the ordering of the columns.
Consequently, we obtain $(RQR)^{k}=R Q^k R = Q^k$, $R e_2 = e_3$ for the $4\times 4$ matrix $Q$ which entails
\begin{equation}
\label{eq:QM}
1^T Q^k e_2  = 1^T R Q^k R e_2 =  1^T Q^k e_3.
\end{equation}
Now, consider
\begin{eqnarray}
F_{k+2} & = & 1^T Q^{k+1} e_1 = 1^T Q^{k-1} Q^2 e_1 = 1^T Q^{k-1} (e_1 + e_3) \nonumber \\
			  & = & 1^T Q^{k-1} (e_3 + e_1) \nonumber \\
				& = & 1^T Q^{k-1} e_2 + 1^T Q^{k-1} e_1 \nonumber \\
				& = & 1^T Q^{k-1} Q e_1 + 1^T Q^{k-1} e_1 \nonumber \\
				& = & F_{k+1} + F_{k}.
\end{eqnarray}

The eigenvalue decomposition of $Q$ yields (details~\cite{Moser2014})
\begin{eqnarray}
\label{eq:xicos}
F_{k} & = &  
\frac{2^{k+2}}{5}\sum_{j=1}^{2}
\cos^k
\left((2j-1) \frac{\pi}{5}\right)
\cos^2\left(\frac{2j-1}{2}\frac{\pi}{5}\right).
\end{eqnarray} 
Due to $\cos(\pi/5) = (1 + \sqrt{5})/4$, $\cos(3 \pi/5) = (1 - \sqrt{5})/4$,
$\cos(\pi/10)^2 = (5 + \sqrt{5})/8$ and 
$\cos(3 \pi/10)^2 = (5 - \sqrt{5})/8$ Equ.~(\ref{eq:xicos}) yields Binet's formula 
\[
F_{k} = F_k = ((1+\sqrt{5})^k - (1-\sqrt{5})^k)/(2^k \sqrt{5}).
\]
By taking advantage of the trigonometric formula
\begin{eqnarray}
\label{eq:tr1}
\cos^m(\theta)
& = & 
\left\{
\begin{array}{lr}
\frac{1}{2^m} 
\left(
\begin{array}{c}
	m\\
      \frac{m}{2}
\end{array}
\right) + 
\frac{2}{2^m}
\sum_{k=0}^{\frac{m}{2}
-1}
 \left(
\begin{array}{c}
	m\\
      k
\end{array}
\right) 
\cos((m-2k) \theta) &  2|m 
 \\ 
\frac{2}{2^m}
\sum_{k=0}^{\frac{m-1}{2}}
 \left(
\begin{array}{c}
	m\\
      k
\end{array}
\right) 
\cos((m-2k) \theta) & 2 \not| m
\end{array}
\right. 
\end{eqnarray}
we are able to replace (after some tedious work) the trigonometric expressions in~(\ref{eq:xicos}) 
in terms of binomial coefficients.
Finally, we obtain for $k$ even
\[
F_{k+1} = 
			\left(
			\begin{array}{c}
				k\\
			\frac{k}{2}
			\end{array}
			\right) +
			2
			\sum_{q=1}^{\left\lfloor \frac{k}{10}\right\rfloor}
 \left(
\begin{array}{c}
	k\\
      \frac{k}{2} - 5 q
\end{array}
\right) 
-
\sum_{q=1, q\, \mbox{\tiny odd}}^{\left\lfloor \frac{k+1}{5}\right\rfloor}
 \left(
\begin{array}{c}
	k+1\\
      \frac{k+1}{2} - \frac{5}{2} q
\end{array}
\right)
\]
and for $k$ odd

\[
F_{k+1} = 
\frac{1}{2}
			\left(
			\begin{array}{c}
				k+1\\
			\frac{k+1}{2}
			\end{array}
			\right) -
2
\sum_{q=1,\mbox{\tiny $q$ odd}}^{\left\lfloor \frac{k}{5}\right\rfloor}
 \left(
\begin{array}{c}
	k\\
\frac{k}{2} - \frac{5}{2} q 
\end{array}
\right) +
\sum_{q=1}^{\left\lfloor \frac{k+1}{10}\right\rfloor}
 \left(
\begin{array}{c}
	k+1\\
      \frac{k+1}{2} - 5 q
\end{array}
\right). 
\]
This summation identity is illustrated in Fig.~\ref{fig:FiboPascal}.

\end{document}